\documentclass[11pt,a4paper]{amsart}
\usepackage[latin1]{inputenc}
\usepackage{amssymb}
\usepackage{amsmath}
\usepackage{latexsym}
\usepackage{epic}
\usepackage{eepic}
\usepackage{graphicx}
\usepackage{psfrag}
\usepackage{longtable}
\usepackage{epsfig}
\usepackage{xypic}
\newtheorem{rem}{Remark}[section]
\newtheorem{prop}{Proposition}[section]

\newtheorem{cor}{Corollary}[section]

\newcommand{\bprf}{{\it Proof.~}}
\newcommand{\eprf}{\hfill $\square$ \bigskip\par}

\newcommand{\pitr}{ \mathbb{P}_3}

\newcommand{\piu}{\mathbb{P}_1}

\newcommand{\ci}{ \mathbb{C}}

\newcommand{\Z}{\mathbb{Z}}

\author{Alessandra Sarti}

\title{Group actions, cyclic coverings and families of K3-surfaces}

\date{\today}

\begin{document}

\maketitle

\begin{abstract}
In this paper we describe six pencils of $K3$-surfaces which have large Picard-number ($15\leq \rho\leq 20$) and contain precisely five singular fibers: four have A-D-E singularities and one is non-reduced. In particular we describe these surfaces as cyclic coverings of the $K3$-surfaces of \cite{basa}. In many cases using this description and lattice-theory we are able to compute the exact Picard-number and to describe explicitly the Picard-lattices.
\end{abstract}



\section{Introduction}
In the last years using different methods (toric geometry, mirror symmetry, etc.) have been constructed and studied many families of $K3$-surfaces with large Picard-Number and small number of singular fibers (see e.g. \cite{dolgachev}, \cite{verryui} and \cite{belcastro}). In these notes using group actions and cyclic coverings we describe some new one. In \cite{basa} the authors describe three pencils of $K3$-surfaces where the generic surface has Picard-number 19 and the pencils contain four singular fibers with singularities of $A-D-E$ type and one non-reduced fiber. The families arise as minimal resolutions of quotients $X_{\lambda}/G$ were $G$ is a subgroup of $SO(4)$ containing the Heisenberg group and $\{X_{\lambda}\}_{\lambda\in\piu}$  is a $G$-invariant pencil of surfaces in $\pitr$, the latter are described in \cite{sa1}. The groups $G$  are the so-called {\it bi-polyhedral groups}: $T\times T$, $O\times O$, $I\times I$, where $T$, $O$ and $I$ denote the tetrahedral group, the octahedral group and the icosahedral group in $SO(3)$ (we recall their precise definition in the first Section), the pencils $X^6_{\lambda}$, $X^8_{\lambda}$, $X^{12}_{\lambda}$ of $T\times T-$, $O\times O-$, resp. $I\times I-$invariant surfaces, have degrees 6, 8, resp. 12. In these notes we consider some special normal subgroup $H$ of $T\times T$ and $O\times O$ ($I\times I$ is simple) so  that the minimal resolutions of the quotients are pencils of $K3$-surfaces, these have then five singular fibers (one is non-reduced) and large Picard-number. In Section 2 we describe all the subgroups $H$ with this property. In the Sections 3, 4, 5 and 6 we describe six new one-dimensional families of $K3$-surfaces: the generic surface of two of the pencils has Picard-number $\rho=19$ and there are four singular surfaces  with $\rho=20$ as in the case of the families $X^6_{\lambda}/T\times T$ or $X^8_{\lambda}/O\times O$ of \cite{basa}. The other four pencils seem to have smaller Picard-number. We give in any case a lower bound for $\rho$ and  in each case except one we could identify surfaces with $\rho=20$. The methods which we use in these Sections are essentially the same as in \cite{basa}. In the Sections 7 and 9, we describe the K3-surfaces as cyclic coverings of the K3-surfaces of \cite{basa}, more precisely let $Y_{\lambda,T\times T}$ and  $Y_{\lambda,O\times O}$ denote the minimal resolutions of the quotients $X^6_{\lambda}/T\times T$ and  $X^8_{\lambda}/O\times O$, then taking two special 3-divisible classes of rational curves in the Neron-Severi group $NS(Y_{\lambda,T\times T})$ and two special 2-divisible classes of rational curves in the Neron-Severi group $NS(Y_{\lambda,O\times O})$ one can construct cyclic coverings (two for each surface) which correspond (up to contract some (-1)-curves) to the minimal resolutions of the surfaces  $X^6_{\lambda}/H$ and  $X^8_{\lambda}/H'$ for some normal subgroup $H$ of $T\times T$ and $H'$ of $O\times O$. By doing cyclic coverings of the latter, one obtains more $K3$-surfaces.
In Section 8 and 9 by using these descriptions and the results of \cite{basa}, Section 6,  we compute  explicitly the Picard-lattice of the $K3$-surfaces of two of the families, more precisely of those where the generic surface has $\rho=19$ and of the four special surfaces with $\rho=20$ contained in these families. We compute also the Picard-lattice of the surfaces  with $\rho=20$ in the other pencils.\\
I would like to thank W. Barth for introducing me to cyclic coverings and for many useful discussions.

\section{Notations and preliminaries}
There are two classical $2:1$ coverings:
\begin{eqnarray*}
\rho:SU(2)\rightarrow SO(3)~~~\hbox{and}~~~\sigma:SU(2)\times SU(2)\rightarrow SO(4).
\end{eqnarray*}
Let $G_i\subset SO(3)$, $i=1,2$ denote the polyhedral group $T$ or $O$. We consider the binary group $\tilde{G_i}:=\rho^{-1}(G_i)\subset SU(2)$ and the $\sigma$-image: $\sigma(\tilde{G_1}\times \tilde{G_2})\subset SO(4)$, which by abuse of notation we denote with $G_1\times G_2$ and we denote an element of $SU(2)\times SU(2)$ and its image in $SO(4)$ by $(p_1,p_2)$. These groups have been studied in \cite{sa1}, there  the group $T\times T$ is called $G_6$ and the group $O\times O$ is called $G_8$.
We will denote by $X^6_{\lambda}=s_6+\lambda q^3$  and by $X^8_{\lambda}=s_8+\lambda q^4$ the pencils of $T\times T$- and of $O\times O$-invariant surfaces in $\pitr$, which are described in \cite{sa1},  
$s_6$ denotes a $T\times T $-invariant homogeneous polynomial of degree six and $s_8$ denotes an $O\times O $-invariant homogeneous polynomial of degree eight, $q:=x_0^2+x_1^2+x_2^2+x_3^2$ is the equation of the quadric $\piu\times \piu$ in $\pitr$. The base locus of the pencils $X^n_{\lambda}$ are $2n$ lines $(n=6,8)$ on the quadric, $n$ in each ruling. These pencils contain exactly four nodal surfaces. We recall the value of the $\lambda_i$ and the number of nodes on $X_{\lambda_i}$ in the table below (cf. \cite{sa1}):

\renewcommand{\arraystretch}{1.3}
\begin{eqnarray*}
\begin{array}{cccc|cccc}
\multicolumn{4}{c|}{n=6} & \multicolumn{4}{c}{n=8}  \\ \hline
\lambda_1 & \lambda_2 & \lambda_3 & \lambda_4 &
\lambda_1 & \lambda_2 & \lambda_3 & \lambda_4   \\
-1 & -\frac{2}{3} & -\frac{7}{12} & -\frac{1}{4} &
-1 & -\frac{3}{4} & -\frac{9}{16} & -\frac{5}{9}  \\
12&48&48&12&24&72&144&96\\
\end{array} 
\end{eqnarray*}
\renewcommand{\arraystretch}{1.0}
We recall also the matrix:
\begin{eqnarray*}
C:=\left(\begin{array}{rrrr}
1&0&0&0\\
0&-1&0&0\\
0&0&-1&0\\
0&0&0&-1
\end{array}
\right)\in O(4),
\end{eqnarray*}
which operates on an element $(p_1,p_2)\in G_1\times G_2$ by:
\begin{eqnarray*}
C^{-1}(p_1,p_2)C=(p_2,p_1).
\end{eqnarray*}
Moreover we specify the following matrices of $SO(4)$:
\begin{eqnarray*}
\begin{array}{ll}
(q_1,1)=
\left( \begin{array} {rrrr}
0 & -1& 0 & 0 \\
1 & 0& 0 & 0 \\
0 & 0& 0 & -1 \\
0 & 0& 1 & 0 
\end{array}\right),&
(1,q_{1})=
\left( \begin{array} {rrrr}
0 & 1& 0 & 0 \\
-1 & 0& 0 & 0 \\
0 & 0& 0 & -1 \\
0 & 0& 1 & 0 
\end{array} \right),
\end{array}
\end{eqnarray*}
\begin{eqnarray*}
\begin{array}{ll}
(q_{2},1)=
\left( \begin{array} {rrrr}
0 & 0& -1 & 0 \\
0 & 0& 0 & 1 \\
1 & 0& 0 & 0 \\
0 & -1& 0 & 0 
\end{array} \right),&
(1,q_{2})=
\left( \begin{array} {rrrr}
0 & 0& 1 & 0 \\
0 & 0& 0 & 1 \\
-1 & 0& 0 & 0 \\
0 & -1& 0 & 0 
\end{array} \right),
\end{array}
\end{eqnarray*}
\begin{eqnarray*}
\begin{array}{ll}
(p_3,1)=\frac{1}{2}
\left( \begin{array} {rrrr}
1& -1& 1& -1 \\
1 & 1& -1 & -1 \\
-1 & 1& 1 & -1 \\
1 & 1& 1 & 1 
\end{array} \right),&
(1,p_3)=\frac{1}{2}
\left( \begin{array} {rrrr}
1 & 1& -1 & 1 \\
-1 & 1& -1 & -1 \\
1 & 1& 1 & -1 \\
-1 & 1& 1 & 1 
\end{array} \right),\\
\end{array}
\end{eqnarray*}

\begin{eqnarray*}
\begin{array}{ll}
(p_4,1)=\frac{1}{\sqrt{2}}
\left( \begin{array} {rrrr}
1& -1& 0& 0 \\
1 & 1& 0 & 0 \\
0 & 0& 1 & -1 \\
0 & 0& 1 & 1 
\end{array} \right),&
(1, p_4)=\frac{1}{\sqrt{2}}
\left( \begin{array} {rrrr}
1 & 1& 0 & 0 \\
-1 & 1& 0& 0 \\
0 & 0& 1 & -1 \\
0 & 0& 1 & 1 
\end{array} \right),
\end{array}
\end{eqnarray*}
Using these matrices the groups have the following generators: 

\begin{center}

\renewcommand{\arraystretch}{1.3}

\begin{eqnarray*}
\begin{array}{c|c}
{\rm Group } &{\rm Generators } \\
\hline
T\times T & (q_2,1), (1,q_2), (p_3,1), (1,p_3)\\
O\times O & (q_2,1), (1,q_2), (p_3,1), (1,p_3), (p_4,1), (1,p_4)\\
\end{array}
\end{eqnarray*}

\end{center}

Denote by $PG$ the image of a subgroup $G\subset SO(4)$ in $\mathbb{P}GL(4)$.
We define the {\it types} of  lines in $\pitr$ which are fixed by elements $(p_1,p_2)\in PG$ of order $2,3$ or $4$ in the following way:
\begin{eqnarray*}
\begin{array}{c|cccc}
order&2&3&4\\
\hline
type&M&N&R\\
\end{array}
\end{eqnarray*}
\section{Normal subgroups}\label{normal}
In \cite{sa2} the author classifies all the subgroups of $SO(4)$ which contain the Heisenberg group $V\times V$. Here we consider all the normal subgroups of $T\times T$ and $O\times O$ which contain the subgroup $V\times V$ and $T\times T$.  The group $I\times I$ does not contain normal subgroups.
We denote by $H$ such a normal subgroup, by $o(H)$ its order and by $i(H)=[G:H]$ the index of $H$ in $G$ ($G=T\times T$ or $O\times O$). We list below all the groups $H$ and their generators, following the notation of \cite{sa2}. Moreover we do not consider separately the groups $H$ and $C^{-1}H C$ or, in general, groups which are conjugate in $O(4)$.

\begin{eqnarray*}
\begin{array}{l|l|l|l|l}
H&\mbox{generators}&o(H)&G&i(H)\\
\hline
T\times V&  (q_1,1),(1,q_1)             &96&T\times T&3\\
          &(p_3,1),(1,p_3)&&&\\
(TT)'&  (q_1,1),(1,q_1)                 &96&T\times T&3\\
&         (q_2,1),(1,q_2)&&&\\        
&          (p_3,p_3)&&&\\   

V\times V& (q_1,1),(1,q_1)             &32&T\times T&9\\
         & (q_2,1),(1,q_2)&&&\\

O\times T&  (q_1,1),(1,q_1)            &576&O\times O&2\\
         & (p_3,1), (1,p_3)&&&\\
         &(p_4,1)&&&\\

(OO)''&  (q_1,1),(1,q_1)               &576&O\times O&2\\
&(p_3,1),(1,p_3)&&&\\
&(p_4q_2,p_4q_2)&&&\\
T\times T& (q_1,1),(1,q_1)             & 288&O\times O&4\\
         &(p_3,1),(1,p_3)&&&\\
\end{array}
\end{eqnarray*}

\section{Fix-points}\label{fix}
We analyze the different kind of fix-points for elements of the subgroups $PH\subset PG$ in the same way as in \cite{basa}. Recall that the elements of the form $(p,1)$ or $(1,p)$ have each two disjoint lines of fix points contained in one ruling, respectively in the other ruling of the quadric (cf. \cite{sa1}, 5.4 p. 439).\\
1) {\it Fix-points on the quadric}. The subgroups $G_1\times 1$ and $1\times  G_2\subset PH$ operate on the two rulings of the quadric and determine orbits of lines, which are also $PH$-orbits of lines. We resume the lengths of the orbits in the following tables. In the first row we write the order of the element which fixes two lines of the orbit:

\begin{eqnarray*}
\begin{array}{l|lll}
{\rm order~~of}~(p,1)&2&3&4\\
\hline
T \times V&6&4,4&-\\
O\times T&12&8&6\\
(TT)'&6&-&-\\
(OO)''&6&8&-\\
V\times V&2,2,2&-&-\\
T\times T&6&4,4&-
\end{array}\; \; \; \; \; \; \; \; \; \;
\begin{array}{l|lll}
{\rm order~~of}~(1,p)&2&3&4\\
\hline
T \times V&2,2,2&-&-\\
O\times T&6&4,4&-\\
(TT)'&6&-&-\\
(OO)''&6&8&-\\
V\times V&2,2,2&-&-\\
T\times T&6&4,4&-
\end{array} 
\end{eqnarray*}

In particular observe that in the case of the groups $(TT)'$ and $(OO)''$ the meeting points of the fix-lines of the two rulings of $\piu\times \piu$ split into three orbits of length 12 and two orbits of length 32, in the other cases these meeting-points form just one orbit.\\
2) {\it Fix-points off the quadric}. We denote by $F_L$ the fix-group of a line L of $\pitr$ in $PH$  and by $H_L$ the stabilizer group of $L$ in $PH$, i.e.
\begin{eqnarray*}
\begin{array}{l}
F_L:=\{h\in PH~s.t.~hx=x~for~all~x\in L\}\\
H_L:=\{h\in PH~s.t.~hL=L\}.\\
\end{array}
\end{eqnarray*}
Moreover denote by $\ell(L)$ the length of the G-orbit of the line $L$ and by 
$g$ a representative of a conjugacy class in $G$:

\begin{eqnarray*}
\begin{array}{c|ccc|cccc|c}
\mbox{group}        & \multicolumn{3}{c|}{T\times V} &\multicolumn{4}{c|}{(TT)'} &
                                     \multicolumn{1}{c}{V\times V} \\ \hline         
g           & (q_1,q_1) & (q_1,q_2) & (q_1,q_3) &(q_1,q_1)&(q_1,q_2)&(q_1,q_3)&(p_3,p_3) & (q_i,q_j)\\
F_L         & \Z_2 & \Z_2 & \Z_2 & \Z_2 & \Z_2&\Z_2&\Z_3 & \Z_2 \\
\mbox{type} & M_1  & M_2  & M_3 & M_1&M_2&M_3&   N & M_{ij} \\
\ell(L)     & 6 & 6 & 6 & 6&6&6 & 16 & 2 \\
|H_L|/|F_L| &  4 &  4 &  4 & 4&4&4  & 1  &  4 \\ 
\end{array}
\end{eqnarray*}
Here we denote by $q_3$ the product of $q_1$ and $q_2$. In the last column of the table the sum runs over $i,j=1,2,3$. In this case we have nine distinct conjugacy classes with just one element.  
\begin{eqnarray*}
\hspace*{-1.5cm}
\begin{array}{c|ccc|cccc|ccc}
\mbox{group}           & \multicolumn{3}{c|}{O\times T} &\multicolumn{4}{c|}{(OO)''} &
                                     \multicolumn{3}{c}{T\times T} \\ \hline         
g           & (q_1,q_1) & (p_3,p_3) & (p_4q_2,q_2) &(p_4,p_4)&(p_3,p_3) & (p_3^2,p_3)& (p_4q_2,p_4q_2)& (q_2,q_2)&(p_3,p_3)& (p_3^2,p_3)\\
F_L         & \Z_2 & \Z_3 & \Z_2 & \Z_4 & \Z_3 & \Z_3&\Z_2&\Z_2&\Z_3&\Z_3 \\
\mbox{type} & M  & N  & M' & R  & N & N'& M  & M  & N  & N'    \\
\ell(L)     & 18 & 32 & 36& 18 & 16 & 16& 72&18&16&16 \\
|H_L|/|F_L| &  8&  3 &  3 & 4  & 8  &  8& 2& 4& 3& 3\\ 
\end{array} 
\end{eqnarray*}

\begin{rem}\label{ss2}
By taking the generator $(p_3^2,p_3)$ for $(TT)'$ instead of $(p_3,p_3)$ we find a group $(TT)''$ which is conjugate in $O(4)$ to $(TT)'$. The singularities in the quotient  $X^6_{\lambda}/(TT)''$ are the same as those on $X^6_{\lambda}/(TT)'$. 
\end{rem}

\section{Quotient singularities}\label{quo}
We consider now the projections:
\begin{eqnarray*}
\begin{array}{cc}
\pi_H:X^6_{\lambda}\longrightarrow X^6_{\lambda}/H, &\pi_{H'}:X^8_{\lambda}\longrightarrow X^8_{\lambda}/H'
\end{array}
\end{eqnarray*}
with $H=T\times V$, $(TT)'$ or $V\times V$; $H'=O\times T$, $(OO)''$ or $T\times T$. In this Section we run the same program as in \cite{basa}, Section 3 and describe the singularities of the quotients (for the details cf. \cite{basa}) .\\
1) {\it Fix-lines on q}. The image in the quotient of the lines of the base locus of the pencils $X^6_{\lambda}$ and $X^8_{\lambda}$ and of the intersection points of the lines of the base locus are smooth. Observe that the points of intersection of the lines of the base locus of the pencils form one orbit under the action of $T\times V$, $V\times V$, $O\times T$ and $T\times T$. In the case of the groups $(TT)'$ we have three orbits and in the case of the group $(OO)''$ we have two orbits, as described in Section \ref{fix}, this means that the lines in the quotient will meet three times and two times. Now we consider the points of intersection of the lines of the base locus with the other fix-lines on $q$. In the table below we do not write the groups $(TT)'$ and $V\times V$ because they do not have other fix-points on $q$ other then the lines of the base locus. We denote by Fix$(P)$ the fix-group in $PG$ of a point $P$. In the next table we write the length and the number of orbits of fix-points, and we describe which kind of singularities do we have in the quotient:  
\begin{eqnarray*}
\begin{array}{c|c|ccc|cc|cc}
\mbox{group}&T\times V& \multicolumn{3}{c|}{O\times T}&\multicolumn{2}{c|}{(OO)''}&\multicolumn{2}{c}{T\times T}\\
\hline
\mbox{Fix}(P)&\Z_3\times\Z_2&\Z_3\times Z_2&\Z_4\times \Z_3&\Z_2\times\Z_3&\Z_2\times\Z_3&\Z_3\times\Z_2&\Z_3\times\Z_2&\Z_2\times\Z_3\\
\mbox{length}&8&48&24&48&48&48&24&24\\
\mbox{number}&6&1&2&2&1&1&2&2\\
\mbox{sing.}&6A_2&1A_1&2A_3&2A_1&1A_1&1A_1&2A_1&2A_1\\
\end{array}
\end{eqnarray*}

2){\it Fix-lines off q}. Denote by  $o(L)$ the order of the fix-group $F_L$ of $L$. The number of points not on $q$ cut out on  $X_{\lambda}$ by   $L$  is: 
\begin{eqnarray*}
\begin{array}{c|c|cc|c|cc|cc|cc}
\mbox{group}         & T\times V & \multicolumn{2}{c|}{(TT)'} &V\times V &\multicolumn{2}{c|}{O\times T}&\multicolumn{2}{c|}{(OO)''}&\multicolumn{2}{c}{T\times T}\\
\hline                                  
o(L)          & 2 &  2 & 3 & 2 & 2 & 3 & 4 &3&2&3\\
\mbox{number} &  4 &4 & 6 & 4 & 8 & 6 &8  &6& 8&6\\
\end{array} 
\end{eqnarray*}
In the next table we show in
each case length and number of $H_L$-orbits, the number and type(s)
of the quotient singularity(ies):
$$ \begin{array}{c|ccc|cccc|c}
\mbox{group}                  & \multicolumn{3}{c|}{T \times V} & \multicolumn{4}{c|}{(TT)'}& V\times V \\
                                                \hline
o(L)                 & 2 & 2 & 2 & 2 & 2& 2 & 3 & 2 \\
\mbox{type}          & M_1 & M_2&M_3& M_1&M_2&M_3& N & M_{ij}\\
\mbox{length}        & 4 & 4 & 4 & 4 & 4 & 4 & 1 & 4\\ 
\mbox{number}        & 1 & 1 & 1 & 1 & 1 & 1 & 6 & 1 \\
\mbox{singularities} & A_1 & A_1 & A_1 & A_1 & A_1 & A_1 & 6A_2 & A_1\\ 
\end{array} $$
$$ \begin{array}{c|ccc|cccc|ccc}
\mbox{group}                  & \multicolumn{3}{c|}{O \times T} & \multicolumn{4}{c|}{(OO)''}& \multicolumn{3}{c}{T \times T}\\
                                                \hline
o(L)                 & 2 & 3 & 2 & 4 & 3& 3 & 2 & 2&3&3 \\
\mbox{type}          & M & N&M'& R&N&N'& M & M&N&N'\\
\mbox{length}        & 8 & 3 & 3 & 4 & 6 & 6 & 2 &4&3&3 \\ 
\mbox{number}        & 1 & 2 & 2 & 2 & 1 & 1 & 4 & 2&2&2 \\
\mbox{singularities} & A_1 & 2A_2 & 2A_1 & 2A_3 & A_2 & A_2 & 4A_1 & 2A_1&2A_2&2A_2\\ 
\end{array} $$
3) {\it The singular surfaces.} We denote by $ns$ the number of nodes on the surfaces and by $F$ the fix-group of a node in $H$. In the table below, we give the number of orbits of nodes and their fix-groups in $PH$, $PH'$ and we describe the singularities in the quotient. For doing this recall the Proposition 3.1 of \cite{basa}:
\begin{prop}
Let $X$ be a nodal surface with $F\subset SO(3)$ the fix-group of the node. Then the image of this node on $X/H$ is a quotient singularity locally isomorphic with $\ci^2/\tilde{F}$, where $\tilde{F}\subset SU(2)$ is the binary group which corresponds to $F$. 
\end{prop}

$$\hspace*{-1.0cm} \begin{array}{c|cccc|cccc|cccc}
\mbox{group}& \multicolumn{4}{c|}{T \times V}& \multicolumn{4}{c|}{(TT)'}& \multicolumn{4}{c}{V\times V}\\
\hline
\lambda&\lambda_1&\lambda_2&\lambda_3&\lambda_4&\lambda_1&\lambda_2&\lambda_3&\lambda_4&\lambda_1&\lambda_2&\lambda_3&\lambda_4\\
ns&12&48&48&12&12&48&48&12&12&48&48&12\\
\mbox{orbit}&1&1&1&1&3&3&1&1&3&3&3&3\\
F&\Z_2\times\Z_2&id&id&\Z_2\times\Z_2&T&\Z_3&id&\Z_2\times\Z_2&\Z_2\times\Z_2&id&id&\Z_2\times\Z_2\\
\hline
\mbox{lines}&1M_1&-&-&1M_1&3M_i&1N&-&1M_1&3M_{ij}&-&-&3M_{ij}\\
\mbox{meeting} &1M_2&&&1M_2&4N&&&1M_2&&&&\\
&1M_3&&&1M_3&&&&1M_3&&&&\\
\hline
\mbox{sing.}&D_4&A_1&A_1&D_4&3E_6&3A_5&A_1&D_4&3D_4&3A_1&3A_1&3D_4\\
\end{array}$$
$$ \hspace*{-1.0cm}\begin{array}{c|cccc|cccc|cccc}
\mbox{group}& \multicolumn{4}{c|}{O\times T}& \multicolumn{4}{c|}{(OO)''}& \multicolumn{4}{c}{T\times T}\\
\hline
\lambda&\lambda_1&\lambda_2&\lambda_3&\lambda_4&\lambda_1&\lambda_2&\lambda_3&\lambda_4&\lambda_1&\lambda_2&\lambda_3&\lambda_4\\
ns&24&72&144&96&24&72&144&96&24&72&144&96\\
\mbox{orbit}&1&1&1&1&2&1&1&2&2&1&1&2\\
F&T&\Z_2\times\Z_2&\Z_2&\Z_3&O&\Z_4&\Z_2&D_3&T&\Z_2&id&\Z_3\\
\hline
\mbox{lines}&3M&1M&1M'&1N&3R&1R&1M&1N(N')&3M&1M&-&1N(N')\\
\mbox{meeting}&4N&2M'&&&4N(N')&&&3M&4N(N')&&&\\
&&&&&6M&&&&&&&\\
\hline
\mbox{sing.}&E_6&D_4&A_3&A_5&2E_7&A_7&A_3&2D_5&2E_6&A_3&A_1&2A_5\\
\end{array}$$
\section{Rational curves}\label{rational}
Let
\begin{eqnarray*}
\mu:Y_{\lambda, H}\longrightarrow X_{\lambda}/H
\end{eqnarray*}
be the minimal resolution of the singularities of $X_{\lambda}/H$, which denotes one of the six quotients of Section \ref{quo}. In the following table we give the number of rational curves coming from the curves of the base locus of $X_{\lambda}$ (denote it by $\nu_1$) and from the resolution of the singularities. The latter are of three kinds: those coming from the intersection points of the lines of the base locus with other fix-lines on $q$, those coming from fix-points which are off $q$ and do not come from nodes of $X_{\lambda}$, and those coming from the nodes. We denote their numbers by $\nu_2$, $\nu_3$ and $\nu_4$, then the total number of rational curves is  $\nu:=\nu_1+\nu_2+\nu_3+\nu_4$. We give the discriminant, $d$,  of the intersection matrix too.\\
{\it 1. The smooth $X_{\lambda}$.}
\begin{eqnarray*}
\begin{array}{c|cccccc}
\mbox{group}&T\times V&(TT)'&V\times V&O\times T&(OO)''&T\times T\\
\hline
\nu_1&4&2&6&3&2&4\\
\nu_2&12&-&-&9&2&4\\
\nu_3&3&15&9&7&14&10\\
\nu&19&17&15&19&18&18\\
d&2^5\cdot 3^3\cdot 5&2^3\cdot 3^6\cdot 5&2^{13}\cdot 5&2^5\cdot 3^3\cdot 7&-2^8\cdot 3^2\cdot 7&-2^2\cdot3^6\cdot 7\\
\end{array}
\end{eqnarray*}
{\it 2. The singular $X_{\lambda}$.} In this case the surfaces $X_{\lambda}$  dont have extra singularities on $q$, hence the number $\nu_1$ and $\nu_2$ remain the same as above and we do not write them again.

$$\hspace*{-1.5cm}\begin{array}{c|cccc|cccc}
\mbox{group}& \multicolumn{4}{c|}{T \times V}& \multicolumn{4}{c}{(TT)'}\\
\hline
\lambda&\lambda_1&\lambda_2&\lambda_3&\lambda_4&\lambda_1&\lambda_2&\lambda_3&\lambda_4\\
\nu_3&-&3&3&-&-&3&15&12\\
\nu_4&4&1&1&4&18&15&1&4\\
\nu&20&20&20&20&20&20&18&18\\
d&-2^4\cdot3^3\cdot5&-2^6\cdot 3^3\cdot 5&-2^6\cdot 3^3\cdot 5&-2^4\cdot3^3\cdot5&-3^3\cdot 5&-2^6\cdot3^3\cdot 5&-2^4\cdot3^6\cdot 5&-2^2\cdot3^6\cdot 5\\
\end{array}$$


$$ \hspace*{-1.0cm}\begin{array}{c|cccc|cccc}
\mbox{group}&\multicolumn{4}{c|}{V\times V}& \multicolumn{4}{c}{O\times T}\\
\hline
\lambda&\lambda_1&\lambda_2&\lambda_3&\lambda_4&\lambda_1&\lambda_2&\lambda_3&\lambda_4\\
\nu_3&-&9&9&-&2&4&5&3\\
\nu_4&12&3&3&12&6&4&3&5\\
\nu&18&18&18&18&20&20&20&20\\
d&-2^{10}\cdot 5&-2^{16}\cdot 5&-2^{16}\cdot 5 &-2^{10}\cdot 5&-2^4\cdot3^2\cdot 7&-2^4\cdot3^3\cdot 7&-2^5\cdot 3^3\cdot 7&-2^6\cdot 3^2\cdot 7\\
\end{array}$$

$$ \begin{array}{c|cccc|cccc}
\mbox{group}& \multicolumn{4}{c|}{(OO)''}& \multicolumn{4}{c}{T\times T}\\
\hline
\lambda&\lambda_1&\lambda_2&\lambda_3&\lambda_4&\lambda_1&\lambda_2&\lambda_3&\lambda_4\\
\nu_3&-&8&12&6&-&8&10&2\\
\nu_4&16&7&3&10&12&3&1&10\\
\nu&20&19&19&20&20&19&19&20\\
d &-2^4\cdot 7&2^7\cdot 3^2\cdot 7&2^8\cdot 3^2\cdot 7&-2^8\cdot 7&-3^4\cdot 7&2^2\cdot 3^6\cdot 7&2^3\cdot 3^6\cdot 7&-2^4\cdot 3^4\cdot 7\\
\end{array}$$

\section {K3-surfaces}\label{k3}
Since the groups $H$ and $H'$ contain the subgroups $V\times V$ of $T\times T$ resp. $T\times T$ of $O\times O$ the projections $\pi_H$ and $\pi_{H'}$ of Section \ref{quo} are ramified on the lines of the base locus of the families $X_{\lambda}$ with ramification index two, and three. A computation as in \cite{basa} Section 5, shows that the quotients are $K3$-surfaces with Picard-number at least the number of rational curves given in Section 5. 
\begin{prop}
The structure of the K3-surfaces $Y_{\lambda, T\times V}$, $Y_{\lambda,O\times T}$ varies with $\lambda$.
\end{prop}
\bprf
\cite{basa} Proposition 5.2 and Lemma 5.2.
\eprf
\begin{cor}
The general $K3$-surface $Y_{\lambda, T\times V}$, $Y_{\lambda,O\times T}$ has Picard-number 19.
\end{cor}
For the other pencils, and in the special cases too, we are not always able to identify exactly the Picard-number, but we can give a lower bound for it, which is the number $\nu$ of rational curves given in the tables above.\\
We are now going to describe another construction for these pencils of $K3$-surfaces.

\section{Cyclic Coverings}\label{cyclic}
Consider now the pairs $G$ and $H$ so that $G/H$ is cyclic, in our cases either $|G/H|=3$ or $|G/H|=2$ (cf. Section \ref{normal}), and consider the map:
\begin{eqnarray*}
\pi: X_{\lambda}/H\longrightarrow X_{\lambda}/G
\end{eqnarray*}
\subsection{The general case} For the moment assume that $X_{\lambda}$ is smooth. The group $G/H$  acts on the points of the fiber $\pi^{-1}(P)$. If the point $P$ is not fixed by $G/H$ then the map is $3:1$ or $2:1$ there. If $P$ is fixed by $G/H$ then we have a singularity on $X_{\lambda}/H$, more precisely an $A_2$ or an $A_1$, now the fiber $\pi^{-1}(P)$ is one point and the map has multiplicity $2$ or $3$ there (cf. \cite{mir} Lemma 3.6 p. 80). Consider now the minimal resolutions $Y_G$ and $Y_H$ of the singular points  of $X_{\lambda}/G$ and $X_{\lambda}/H$. As shown in \cite{basa} and in Section \ref{k3} these  surfaces are K3 and we have a rational map:
\begin{eqnarray*}
\gamma:Y_H---\rightarrow Y_G
\end{eqnarray*}
which is $3:1$ or $2:1$. Observe that this map is not defined over the $(-2)$-rational curves in the blow up of the singular points of $X_{\lambda}/G$ which comes from fix-points of $G/H$ on $X_{\lambda}/H$.\\
In this Section we describe the map $\gamma$ in another way, more precisely by using cyclic coverings. For the general theory about 2-cyclic coverings and 3-cyclic coverings we send back to the article \cite{nikulin} of Nikulin and to the articles \cite{mir1} of Miranda and \cite{tan} of Tan. For the convenience of the reader in Figure 1. we recall the configurations of (-2)-rational curves on the smooth surfaces $Y_{T\times T}$ and on $Y_{O\times O}$ given in \cite{basa}.
\begin{figure}[t]

\input{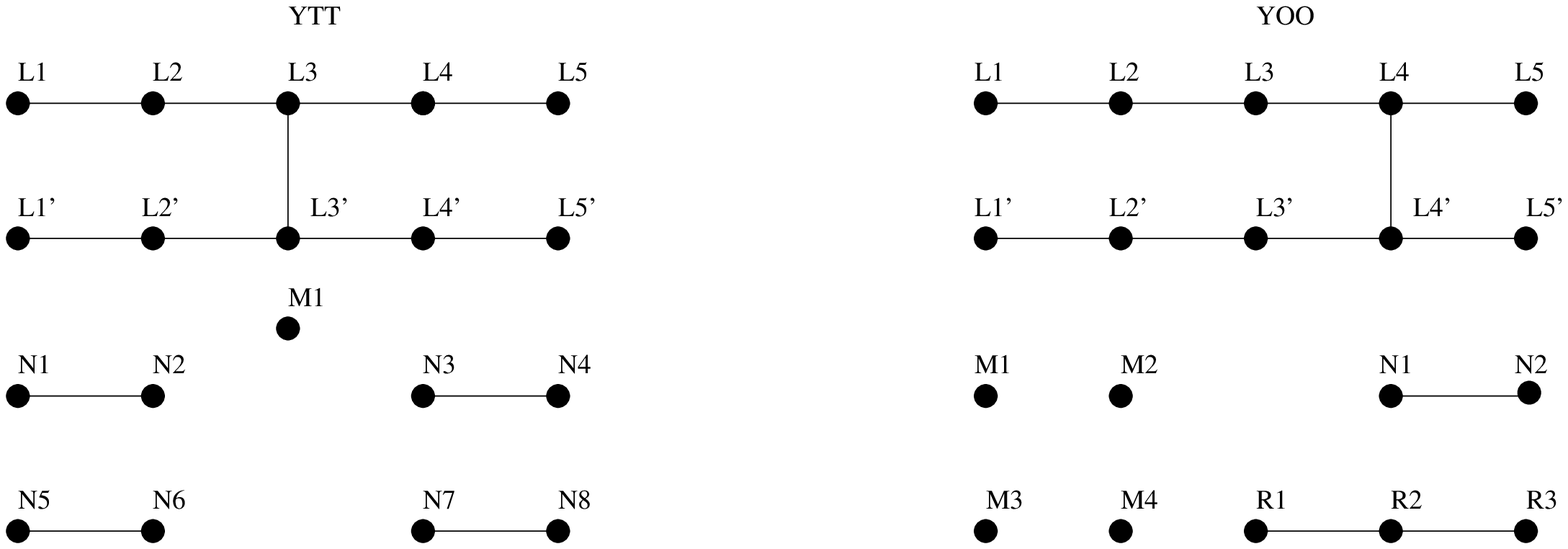}

\begin{center}
\small{Fig. 1}
\end{center}

\end{figure}

\noindent
By \cite{basa}, Proposition 6.1, we have the following 3-divisible classes in  $NS(Y_{T\times T})$:
\begin{eqnarray*}
\mathcal{L}:=L_1-L_2+L_4-L_5+N_1-N_2+N_3-N_4+N_5-N_6+N_7-N_8,\\
\mathcal{L'}:=L_1'-L_2'+L_4'-L_5'+N_1-N_2+N_3-N_4-N_5+N_6-N_7+N_8,
\end{eqnarray*}
and also:
\begin{eqnarray*}
\mathcal{L}-\mathcal{L'}=L_1-L_2+L_4-L_5-L_1'+L_2'-L_4'+L_5'+2(N_5-N_6+N_7-N_8),\\
\mathcal{L}+\mathcal{L'}=L_1-L_2+L_4-L_5+L_1'-L_2'+L_4'-L_5'+2(N_1-N_2+N_3-N_4).
\end{eqnarray*}
Making reduction modulo three we find the classes:
\begin{eqnarray*}
\mathcal{M}:=L_1-L_2+L_4-L_5-L_1'+L_2'-L_4'+L_5'-(N_5-N_6+N_7-N_8),\\
\mathcal{M'}:=L_1-L_2+L_4-L_5+L_1'-L_2'+L_4'-L_5'-(N_1-N_2+N_3-N_4).
\end{eqnarray*}
In $NS(Y_{O\times O})$ we have the following $2$-divisible classes:
\begin{eqnarray*}
\mathcal{L}:=L_1+L_3+L_5+M_1+M_3+M_4+R_1+R_3,\\
\mathcal{L'}:=L_1'+L_3'+L_5'+M_2+M_3+M_4+R_1+R_3.
\end{eqnarray*}
Consider also the classes $\mathcal{L}+\mathcal{L'}$ and $\mathcal{L}-\mathcal{L'}$, which after reduction modulo $2$ are the same as:
\begin{eqnarray*} 
\mathcal{M}:=L_1+L_3+L_5+L_1'+L_3'+L_5'+M_1+M_2.
\end{eqnarray*}
These classes consist of six disjoint $A_2$-configurations of curves and of eight disjoint $A_1$-configurations of curves (in accord to \cite{tan} and \cite{nikulin}). These are the resolutions of $A_2$ and $A_1$ singularities of $X_{\lambda}/G$ which arise by doing the quotient of $X_{\lambda}/H$ by $G/H$. We construct the 3-cyclic coverings and the 2-cyclic coverings by using the divisors $\mathcal{L},\mathcal{L'},
\mathcal{M},\mathcal{M'}$. To avoid to produce singularities by doing this, we first blow up the meeting points of the $A_2$-configurations.  Call $Y^{0}_{T\times T}$ the surface which we obtain after these blow-ups. The meeting points are replaced by $(-1)$-curves and the two $(-2)$-curves become now $(-3)$-curves. Denote by $\phi: Y^{1}_{T\times T}\longrightarrow  Y^{0}_{T\times T}$ the 3-cyclic covering with branching divisor $\mathcal{L}$, $\mathcal{L}'$ or $\mathcal{M}$, $\mathcal{M}'$:

\begin{prop}\label{configu}
A configuration of curves on $Y^0_{T\times T}$:



   \begin{center}
\vspace*{0.5cm}
   \begin{psfrags}
     \psfrag{a}{$-3$}
     \psfrag{b}{$-1$}
     \psfrag{L1}{$L_1$}
     \psfrag{L2}{$L_2$}
     \psfrag{L3}{$L_3$}
     \psfrag{L4}{$L_4$}
     \psfrag{L5}{$L_5$}
\psfrag{L1'}{$L_1'$}
     \psfrag{L2'}{$L_2'$}
     \psfrag{L3'}{$L_3'$}
     \psfrag{L4'}{$L_4'$}
     \psfrag{L5'}{$L_5'$}
     \psfrag{L2''}{$L_2''$}
\psfrag{M1}{$M_1$}
\psfrag{M2}{$M_2$}
\psfrag{M3}{$M_3$}
\psfrag{M4}{$M_4$}
 \psfrag{L0}{$L_0$}
\psfrag{N1}{$N_1$}
\psfrag{N2}{$N_2$}
\psfrag{N3}{$N_3$}
\psfrag{N4}{$N_4$}

\psfrag{R1}{$R_1$}
\psfrag{R2}{$R_2$}
\psfrag{R3}{$R_3$}
\psfrag{R1'}{$R_1'$}
\psfrag{R2'}{$R_2'$}
\psfrag{R3'}{$R_3'$}

     \includegraphics[width=4cm]{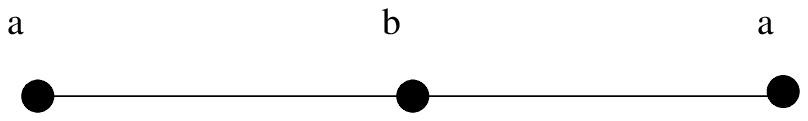}
   \end{psfrags}
    \vspace*{0.5cm}
\end{center}

becomes a configuration: 



   \begin{center}
   \vspace*{0.5cm}
   \begin{psfrags}
     \psfrag{a}{$-3$}
     \psfrag{b}{$-1$}
     \psfrag{L1}{$L_1$}
     \psfrag{L2}{$L_2$}
     \psfrag{L3}{$L_3$}
     \psfrag{L4}{$L_4$}
     \psfrag{L5}{$L_5$}
\psfrag{L1'}{$L_1'$}
     \psfrag{L2'}{$L_2'$}
     \psfrag{L3'}{$L_3'$}
     \psfrag{L4'}{$L_4'$}
     \psfrag{L5'}{$L_5'$}
     \psfrag{L2''}{$L_2''$}
\psfrag{M1}{$M_1$}
\psfrag{M2}{$M_2$}
\psfrag{M3}{$M_3$}
\psfrag{M4}{$M_4$}
 \psfrag{L0}{$L_0$}
\psfrag{N1}{$N_1$}
\psfrag{N2}{$N_2$}
\psfrag{N3}{$N_3$}
\psfrag{N4}{$N_4$}

\psfrag{R1}{$R_1$}
\psfrag{R2}{$R_2$}
\psfrag{R3}{$R_3$}
\psfrag{R1'}{$R_1'$}
\psfrag{R2'}{$R_2'$}
\psfrag{R3'}{$R_3'$}

     \includegraphics[width=4cm]{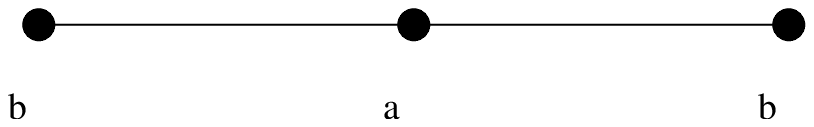}
   \end{psfrags}
\vspace*{0.5cm}
\end{center}

on  $Y^{1}_{T\times T}$.
\end{prop}
\bprf
We do the computation for one configuration of curves  $L_1-L_2$,  this is the same in the other cases. Denote again by $L_1$ and $L_2$ the curves on $Y^0_{T\times T}$ which now are $(-3)$-curves and denote by $E$ the exceptional $(-1)$-curve. By the properties of cyclic coverings we have $\phi^*L_i=3\tilde{L_i}$. Where $\tilde{L_i}$ is the strict transform of $L_i$. Then: 
\begin{eqnarray*}
9(\tilde{L_i})^2=(\phi^* L_i)^2=(\mbox{deg}\phi) L_i^2=-9.
\end{eqnarray*} 
Hence $(\tilde{L_i})^2=-1$. Since $E\cdot (L_1-L_2)=0$ the map $\phi$ is not ramified on $E$ and the restriction $\phi_{|\tilde{E}}$ is $3:1$ onto $E$. Hence we have  $\phi^* E=\tilde{E}$ and $\tilde{E}^2=(\phi^* E)^2=(\mbox{deg}\phi) E^2=3 E^2=-3$.
\eprf
Our surface $Y^{1}_{T\times T}$ is now no more minimal. By blowing down the $(-1)$-curves, the curve $\tilde{E}$ becomes also a $(-1)$-curve so we blow it down too. We call the new smooth surfaces $T_\mathcal{L}$, $T_{\mathcal{L}'}$ resp. $T_\mathcal{M}$, $T_{\mathcal{M}'}$.\\
On $Y_{O\times O}$ we construct directly the 2-cyclic coverings with branching divisors $\mathcal{L}$ and $\mathcal{M}$, because these consist of disjoint smooth rational curves. A computation as in Proposition \ref{configu} shows that on the coverings the $(-2)$-curves of the branching divisors became $(-1)$-curves, hence we can blow  them down. We call the surfaces which we obtain $O_{\mathcal{L}}$, $O_{\mathcal{L}'}$ and $O_{\mathcal{M}}$.\\ 
By Proposition \ref{TO} below and Remark \ref{altrigruppi} follows that the surfaces $T_\mathcal{L}$, $T_{\mathcal{L}'}$, $T_\mathcal{M}$ and $O_{\mathcal{L}}$,  $O_{\mathcal{L}'}$, $O_{\mathcal{M}}$ are minimal K3-surfaces and by construction are exactly the surfaces $Y_{H}$ which are obtained as the minimal resolutions of $X_{\lambda}/H$. We have a commutative diagram as in Figure 2.
\begin{figure}[t]

$$ \begin{array}{ccccc}
&Y^1_G&\stackrel{\phi}{\longrightarrow}&Y^0_G&\\
\swarrow&&&&\searrow\\
Y_H&&-\stackrel{\gamma}{-}\rightarrow&&Y_G\\
\searrow&&&&\swarrow\\
&X_\lambda/H&\stackrel{\pi}{\longrightarrow}&X_\lambda/G&\\
\end{array} $$

\begin{center}
\small{Fig. 2}
\end{center}

\end{figure}

\noindent
In the diagram in the case of $G=O\times O$, we have $Y^0_G=Y_G$ and the map to $Y_G$ is the identity. This construction allow us to describe in another way the $(-2)$-curves on $Y_{H}$  which in some cases determine the rank of the Neron-Severi group of the $K3$-surfaces.
\begin{prop}\label{TO}
The surfaces $T_\mathcal{L}$, $T_\mathcal{M}$, $O_\mathcal{L}$, $O_\mathcal{M}$ are $K3$, containing the following numbers of independent rational $(-2)$ curves coming from the rational curves of  $Y_{T\times T}$ and $Y_{O\times O}$ of Figure 1:

\begin{eqnarray*}
\begin{array}{r|r|r|r}
T_\mathcal{L}& T_\mathcal{M} &O_\mathcal{L} & O_\mathcal{M}\\
\hline
19&17&19&18
\end{array}
\end{eqnarray*}
Moreover we have the following isomorphism of surfaces:

\begin{eqnarray*}
\begin{array}{r|r|r|r}
T_\mathcal{L}& T_\mathcal{M} &O_\mathcal{L} & O_\mathcal{M}\\
\hline
 Y_{T\times V} &Y_{ (TT)'}&Y_{O\times T} &Y_{(OO)''}
\end{array}
\end{eqnarray*}
We have these isomorphisms up to conjugate with $C$ the group $T\times V$, up to change $(p_3,p_3)$ with $(p_3^2,p_3)$ in $(TT)'$ and up to conjugate with $C$ the group $O\times T$. 

\end{prop}
\bprf
1. Observe that the curves $L_1,L_2,L_4,L_5,N_1,N_2,N_3,N_4,N_5,N_6,N_7,N_8$ ``disappear'' on $T_\mathcal{L}$, as the curves $L_1,L_2,L_4,L_5,L_1',L_2',L_4',L_5', N_5,N_6,N_7,N_8$ ``disappear'' on $T_\mathcal{M}$. Moreover since $L_3\cdot \mathcal{L}=0$ the map $\phi$ is not ramified on $L_3$. In the same way since $L_3\cdot\mathcal{M}=L_3'\cdot \mathcal{M}=0$ the map $\phi$ is not ramified on $L_3$, $L_3'$ and the two rational curves meet in three distinct points. Observe that the $(-2)$-curves of $Y'_{T\times T}$ remain $(-2)$-curves on $T_{\mathcal{L}}$ and on $T_{\mathcal{M}}$ and we have the configurations of $(-2)$-rational curves given in Figure 3.
\begin{figure}[t]

\input{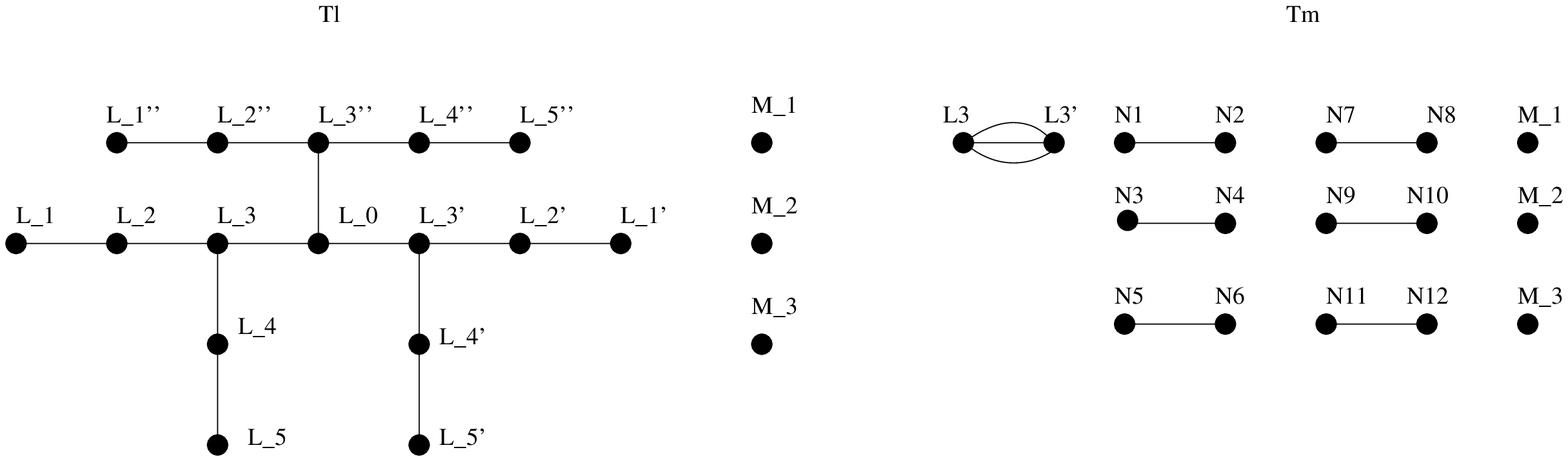}

\begin{center}
\small{Fig. 3}
\end{center}

\end{figure}

\noindent
These are $19$ and $17$ curves. The discriminants are the following:
\begin{eqnarray*}
\begin{array}{c|ccc|c}
               &d(L)                &d(M)      &d(N)  &d\\
\hline
T_{\mathcal{L}}&-2^2\cdot 3^3\cdot 5&-2^3      &      &2^5\cdot 3^3\cdot 5\\
T_{\mathcal{M}}&-5                  &-2^3      &3^6   &2^3\cdot 3^6\cdot 5\\
\end{array}
\end{eqnarray*}
Since it is not zero these classes are independent in the Neron-Severi group.\\
In the same way as before the curves $L_1,L_3,L_5,M_1,M_3,M_4,R_1,R_3$ disappear on $O_{\mathcal{L}}$ as the curves $L_1, L_3, L_5,L_1',L_3',L_5',M_1,M_2$ disappear on $O_{\mathcal{M}}$. Since $L_2\cdot \mathcal{L}=L_4\cdot\mathcal{L}=R_2\cdot \mathcal{L}=0$ and  $L_2\cdot\mathcal{M}=L_4\cdot\mathcal{M}=L_2'\cdot\mathcal{M}=L_4'\cdot\mathcal{M}=0$ the 2-cyclic coverings are not ramified there, and $L_4,L_4'$ meet in two distinct points. Here the configurations of $(-2)$-rational curves are as in Figure 4.\\
These are 19 and 18 curves. The discriminants are the following:
\begin{eqnarray*}
\begin{array}{c|cccc|c}
               &d(L)                &d(M)      &d(N)  &d(R)     &d\\
\hline
O_{\mathcal{L}}&-2^2\cdot 3\cdot 7&2^2         & 3^2  &-2       & 2^5\cdot 3^3\cdot 7 \\
O_{\mathcal{M}}&-7                 &2^4      &3^2     &2^4      &-2^8\cdot 3^2\cdot 7\\
\end{array}
\end{eqnarray*}
These are not zero hence the previous classes are independent in the Neron-Severi group.\\
Finally observe that the canonical bundles on the surfaces  are trivial, hence the surfaces are abelian or $K3$. Since in each case the rank of the Neron-Severi group is at least 17, it turns out that these are $K3$-surfaces and by construction are isomorphic to the surfaces announced in the statement.
\eprf
\begin{figure}[t]

\input{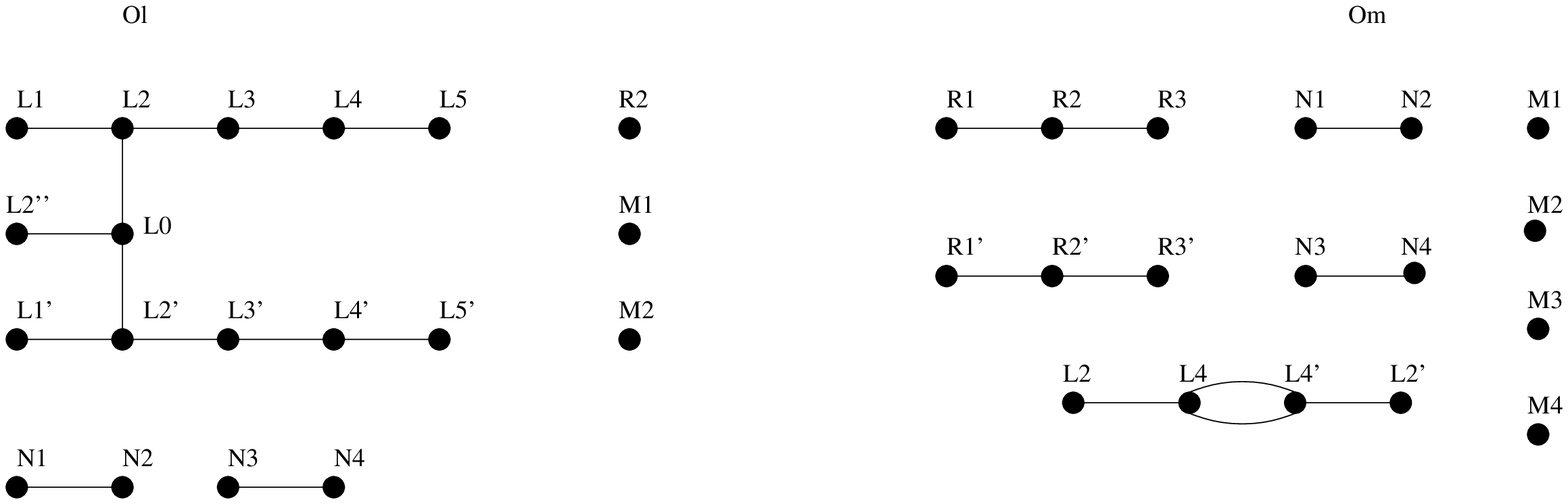}

\begin{center}
\small{Fig. 4}
\end{center}

\end{figure}

\noindent
\begin{rem}\label{altrigruppi}
1. By taking the coverings $T_\mathcal{L'}$ resp. $T_\mathcal{M'}$ one gets  ``almost'' the same surfaces as before, more precisely one has to conjugate the normal subgroup $T\times V$ with $C$, and to interchange $(p_3,p_3)$ in $(TT)'$ with $(p_3^2,p_3)$, getting the group $(TT)''$ (cf. Remark \ref{ss2}).\\
2. By taking the covering $O_{\mathcal{L}'}$ one gets a K3-surface which is isomorphic with $Y_{ T\times O}$.\\
We do not discuss these surfaces separately.\\
\end{rem}

\subsection{The special cases}\label{special}
 We use the same notation of \cite{basa}. We denote by $X_{6,1}$, $X_{6,2}$, $X_{6,3}$, $X_{6,4}$ the four singular surfaces in the pencil $X^6_{\lambda}$ and by $X_{8,1}$, $X_{8,2}$, $X_{8,3}$, $X_{8,4}$ the four singular surfaces in the pencil $X^8_{\lambda}$. For the convenience of the reader we recall in Figure 5 also the graphs of the rational curves in the resolutions of the quotient, which come from singular points outside the quadric. In Figure 5  we do not draw separately the graphs of these  curves for the surfaces $Y^{(6,1)}_{T\times T}$ and $Y^{(6,4)}_{T\times T}$ resp. of $Y^{(6,2)}_{T\times T}$  and $Y^{(6,3)}_{T\times T}$ since these looks equal, however one has to replace the curves $N_1,~N_2,~N_3,~N_4$ by the curves $N_5,~N_6,~N_7,~N_8$.
\begin{figure}[t]

\input{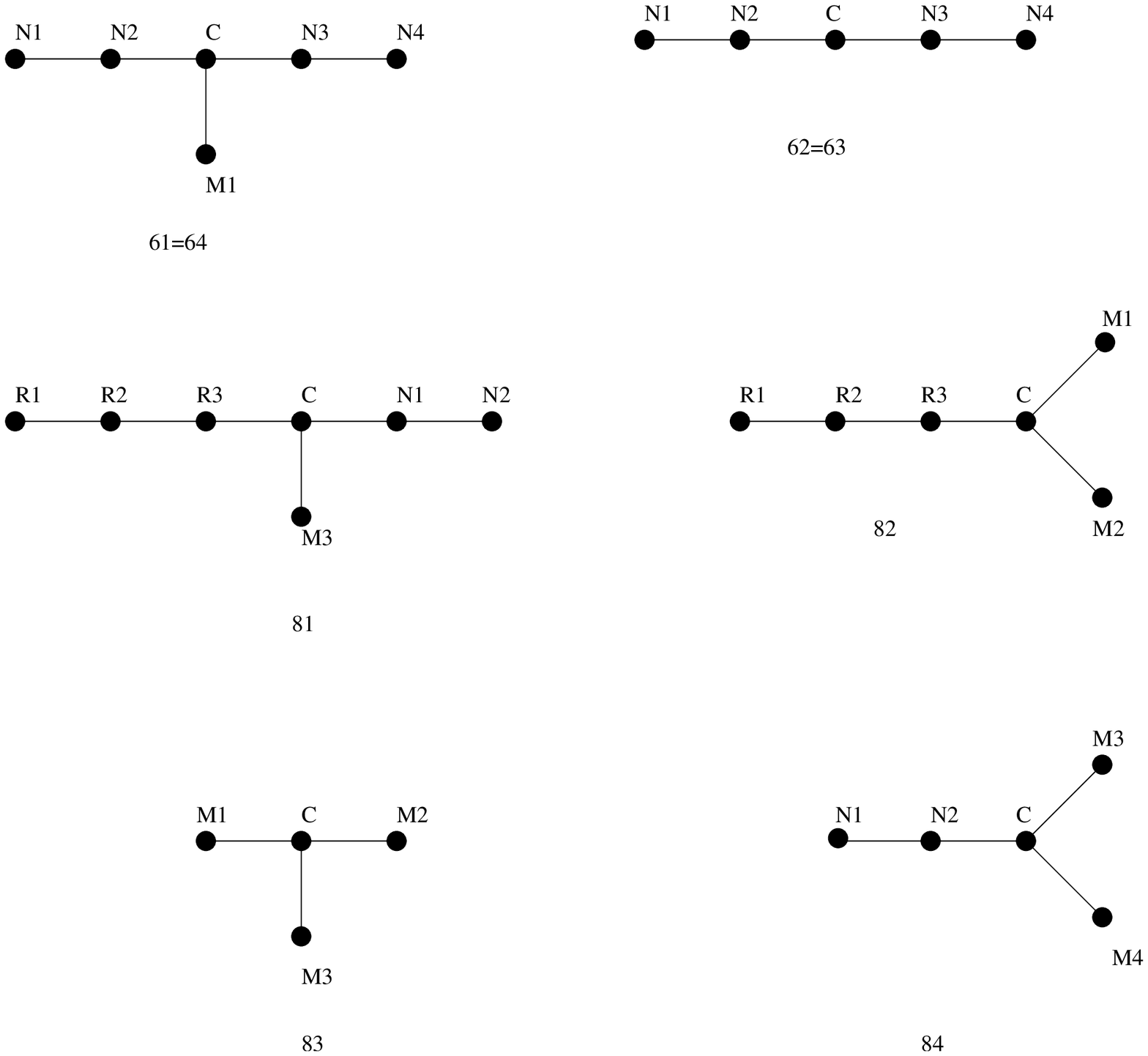}

\begin{center}
\small{Fig. 5}
\end{center}

\end{figure}

\noindent
Now consider the projection $\pi$ of Section \ref{cyclic}. In these cases, the situation is a little more complicated. Now in the counterimage $\pi^{-1}(P)$ of some singular point $P$ of $X_{\lambda}/G$  coming from the $A_1$-singularities of $X_{\lambda}$ we have singularities on $X_{\lambda}/H$ too. In the following table we give the singularities in the quotient $X_{\lambda}/G$, $G=T\times T$ or $O\times O$ and the type and the number of singularities in the counterimage on $X_{\lambda}/H$:\\
\begin{eqnarray*}
\begin{array}{c|cccc|c|cccc}
         &6,1&6,2&6,3&6,4&               &8,1&8,2&8,3&8,4\\
\hline
T\times T&E_6&A_5&A_5&E_6&     O\times O& E_7&D_6&D_4&D_5\\
\hline
T\times V&D_4&A_1&A_1&D_4&     O\times T&E_6&D_4&A_3&A_5\\
(TT)'    &3E_6&3A_5&A_1&D_4&    (OO)''  &2E_7&A_7&A_3&2D_5\\
V\times V&3D_4&3A_1&3A_1&3D_4& T\times T&2E_6&A_3&A_1&2A_5\\
\end{array}
\end{eqnarray*}

By resolving the quotients we get again a map like $\gamma$ in Section \ref{cyclic}. We can describe this map as there by using cyclic coverings. We distinguish two cases:\\
{\it 1. The case of $T\times T$.}
One constructs the  3-cyclic covering as in the general case by using the divisors $\mathcal{L}$, $\mathcal{L'}$, $\mathcal{M}$, $\mathcal{M'}$, these are in the case of the singular surfaces $3$-divisible too. Then one blows down the $(-1)$-curves. The graphs of the $(-2)$-rational curves on the special surfaces and not coming from the lines of the base locus of $X_{\lambda}^6$ (these are the same as in the Figure 3) looks as in Figure 6 (cf. also \cite{basa}, Section 4.2). Observe that the graphs of the curves on $T_{\mathcal{L}}^{(6,1)}$ and on $T_{\mathcal{L}}^{(6,4)}$, resp. on $T_{\mathcal{L}}^{(6,2)}$ and on $T_{\mathcal{L}}^{(6,3)}$ are the same.
\begin{figure}[t]

\input{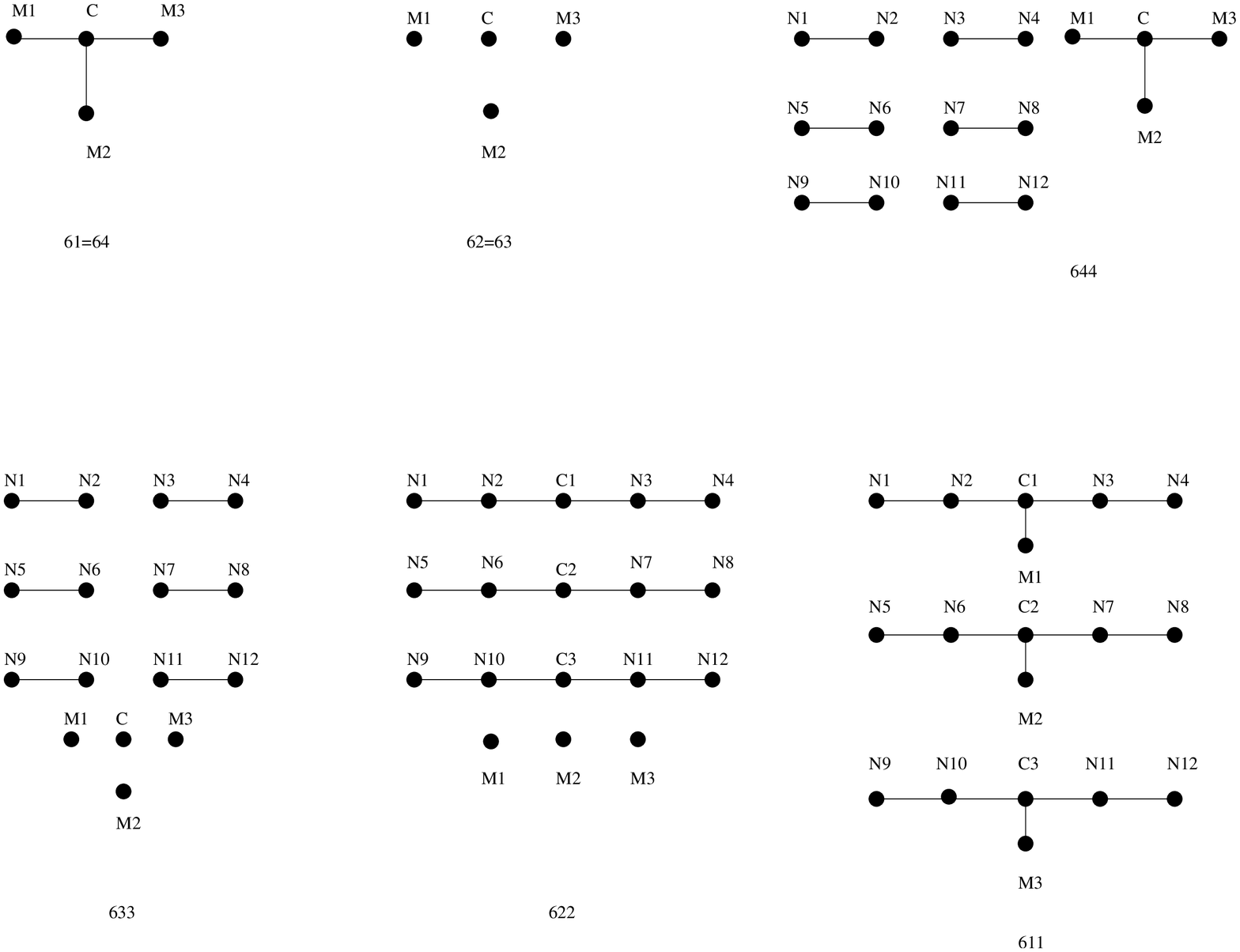}

\begin{center}
\small{Fig. 6}
\end{center}

\end{figure}

\noindent
This together with the $(-2)$-curves coming from the base locus of the pencil $X_{\lambda}^6$ gives the following numbers of independent curves in the Neron-Severi group:
\begin{eqnarray*}
\begin{array}{r|r|r|r|r}
\mathcal{L}&6,1&6,2&6,3&6,4\\
\hline
&20&20&20&20
\end{array},~~~~~~~
\begin{array}{r|r|r|r|r}
\mathcal{M}&6,1&6,2&6,3&6,4\\
\hline
&20&20&18&18
\end{array}
\end{eqnarray*}
Observe that by taking the divisor $\mathcal{L'}$ one obtains the same graphs, by taking the divisor $\mathcal{M'}$ one finds the following isomorphisms of graphs and, by construction, of surfaces:
\begin{eqnarray*}
T^{(6,1)}_{\mathcal{M'}}\cong T^{(6,4)}_{\mathcal{M}},~~~T^{(6,2)}_{\mathcal{M'}}\cong T^{(6,3)}_{\mathcal{M}},~~~T^{(6,3)}_{\mathcal{M'}}\cong T^{(6,2)}_{\mathcal{M}},~~~T^{(6,4)}_{\mathcal{M'}}\cong T^{(6,1)}_{\mathcal{M}}.
\end{eqnarray*}
In fact one sees that by taking $\mathcal{L}$ or $\mathcal{L'}$ one obtains the special $K3$-surfaces in the family $Y_{\lambda, T\times V}$ (up to conjugate the group $T\times V$ with the matrix $C$), by taking $\mathcal{M}$ or $\mathcal{M'}$ one obtains the special $K3$-surfaces in the covering $Y_{\lambda,(TT)'}$ or $Y_{\lambda,(TT)''}$. \\
{\it 2. The case of $O\times O$.}
We take now the divisors $\mathcal{L}$, $\mathcal{L'}$, $\mathcal{M}$ and do $2$-cyclic coverings, then blow down the $(-1)$-curves. The graphs of $(-2)$-rational curves  on the special surfaces  not coming from the lines of the base locus of the pencil $X_{\lambda}^8$ (these looks like in Figure 4) are as in Figure 7.
\begin{figure}[t]

\input{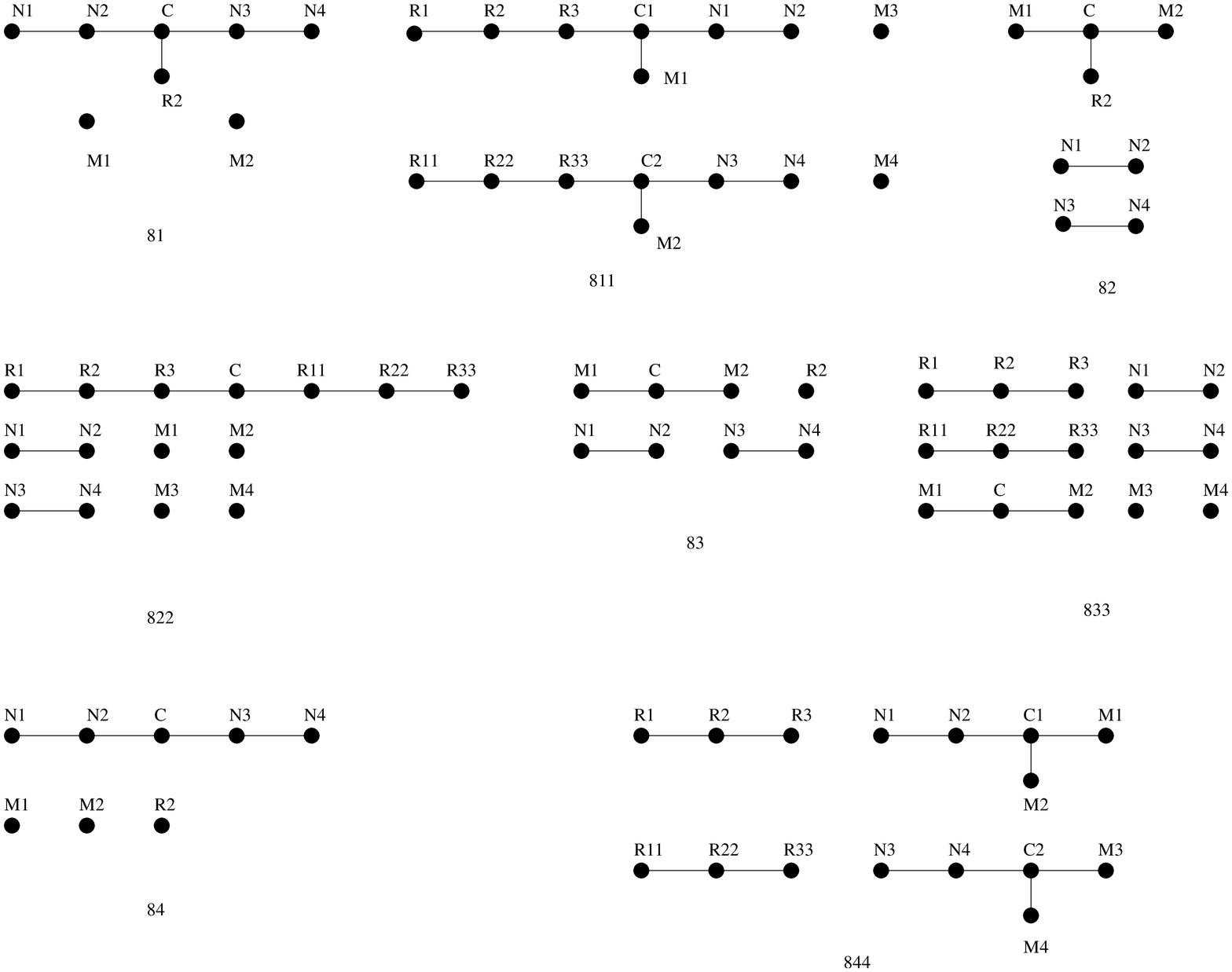}

\begin{center}
\small{Fig. 7}
\end{center}

\end{figure}

\noindent
These together with the $(-2)$-curves coming from the base locus of the pencil $X^8_{\lambda}$ give the following numbers of independent curves in the Neron-Severi group:
\begin{eqnarray*}
\begin{array}{r|r|r|r|r}
\mathcal{L}&8,1&8,2&8,3&8,4\\
\hline
&20&20&20&20
\end{array},~~~~~~~
\begin{array}{r|r|r|r|r}
\mathcal{M}&8,1&8,2&8,3&8,4\\
\hline
&20&19&19&20
\end{array}
\end{eqnarray*}
Similarly to the general case, by taking  $\mathcal{L}$ or $\mathcal{L'}$ we find the singular surfaces in the family $Y_{\lambda,O\times T}$ (up to conjugate $O\times T$ with the matrix $C$) and by taking $\mathcal{M}$ we find the singular surfaces in the family $Y_{\lambda,(OO)''}$.\\
In both the cases 1 and 2 the discriminants of the lattices spanned by the previous curves are those given in Section \ref{rational}.

\section{Picard-lattices}

We compute the Picard-lattices of the general K3-surface in the families  $T_{\mathcal{L}}$,  $O_{\mathcal{L}}$ and of the special surfaces with $\rho=20$ in each pencil. First we recall some facts. Denote by $W$ the lattice spanned by the curves of Proposition \ref{TO}. If $W$ is not the total Picard-lattice, which we call $NS$ there is an integral lattice $W'$ s.t. $W\subset W'\subset NS$ with $p:=[W':W]$ a prime number. Denote by $d(W)$, $d(W')$ the discriminant of the lattices $W,W'$. Since $[W':W]^2=d(W)\cdot d(W')^{-1}$ (cf. \cite{bpv} Lemma 2.1, p. 12)  we find that $p^2$ divides the discriminant of $W$. Denote by $(W^{\vee}/W)^p$ the $p$-subgroup of $(W^{\vee}/W)$$\subset$$NS^{\vee}/NS$ and denote by T the transcendental lattice orthogonal to the Picard-lattice. Since  the discriminant groups $T^{\vee}/T$ and $NS^{\vee}/NS$ are isomorphic (cf. e.g. \cite{bpv} p. 13 Lemma 2.5), they have the same rank which is $\leq$ rk$(T)$. It follows that also rk$(W^{\vee}/W)^p\leq$ rk$(T)$.\\

\subsection{The general case}
\begin{prop}\label{picard}
1. The class:
\begin{eqnarray*}
\bar{L'}:=L_1-L_2+L_4-L_5+L_1'-L_2'+L_4'-L_5'+L_1''-L_2''+L_4''-L_5''\\
\end{eqnarray*}
is 3-divisible in $NS(T_{\mathcal{L}})$ and the classes:
\begin{eqnarray*}
h_1:=L_1+L_3+L_5+L_1'+L_3'+L_5'+M_1+M_2, \\
h_2:=L_1+L_3+L_5+L_1''+L_3''+L_5''+M_1+M_3 \\
\end{eqnarray*}
are 2-divisible in $NS(T_{\mathcal{L}})$. These together with the $19$ curves of Proposition \ref{TO} span a 19-dimensional lattice with discriminant $2\cdot 3\cdot 5$.\\
2. The class:
\begin{eqnarray*}
\bar{L'}:=L_1+L_3+L_5+L_1'+L_3'+L_5'+M_1+M_2
\end{eqnarray*}
is 2-divisible in $NS(O_{\mathcal{L}})$ and the class:
\begin{eqnarray*}
k_1:=L_1-L_2+L_4-L_5-L_1'+L_2'-L_4'+L_5'+N_1-N_2+N_3-N_4\\
\end{eqnarray*}
is 3-divisible. These together with the 19 curves of Proposition \ref{TO} span a 19-dimensional lattice with discriminant $2^3\cdot 3\cdot 7$.\\
\end{prop}
\bprf
1. The discriminant of the lattice generated by the 19 curves is $2^5\cdot 3^3\cdot 5$ hence we can have 2-divisible classes or 3-divisible classes. The divisor $\bar{L'}$ is 3-divisible since it is the pull back of the divisor $\mathcal{L}'$ on $Y_{T\times T}$ which is 3-divisible too. And we cannot have more 3-divisible classes. If there are no 2-divisible classes then the group $(W^{\vee}/W)^2$ would contain the classes $M_1/2$, $M_2/2$, $M_3/2$, $(L_1+L_3+L_5+L_1'+L_3'+L_5')/2$, $(L_1+L_3+L_5+L_1''+L_3''+L_5'')/2$, $(L_1'+L_3'+L_5'+L_1''+L_3''+L_5'')/2$ which are  independent classes with respect the intersection form. Since the rank of $(W^{\vee}/W)^2$ is less or equal as the rank of $T^{\vee}/T$ which is at most three, it can not happen  that we find five classes as before. Hence some combination of them must be  contained in the Neron-Severi group. So we have
\begin{eqnarray*}
\frac{1}{2}(\lambda(L_1+L_3+L_5)+\lambda'(L_1'+L_3'+L_5')+\lambda''(L_1''+L_2''+L_3'')+\mu_1 M_1+\mu_2 M_2+\mu_3 M_3)\in NS
\end{eqnarray*}
for some parameters $\lambda, \lambda',\lambda'',\mu_1,\mu_2, \mu_3\in \Z_2$.\\
By Nikulin \cite{nikulin} such a 2-divisible set contains 8 curves. So putting $\lambda''=0$ and $\mu_3=0$ we get the divisor $h_1$, putting $\lambda'=0$ and $\mu_2=0$ we get the divisor $h_2$. The discriminant of the lattice $W$ together with these three classes now change into $2\cdot 3\cdot 5$, hence we cannot have more torsion classes.\\
2. Again the class $\bar{L'}$ is the pull back of the class $\mathcal{M}'$ on $Y_{O\times O}$ hence  it is 2-divisible. If there are no 3-divisible classes then the group $(W^{\vee}/W)^3$ would contain the classes $N_1-N_2/3$, $N_3-N_4/3$ and $(L_1-L_2+L_4-L_5+L_1'-L_2'+L_4'-L_5')/3$ which are independent. We have these three independent classes also on the special surfaces $O_\mathcal{L}^{(8,2)}$ and on $O_\mathcal{L}^{(8,3)}$. Hence we would have also there rk$(W^{\vee}/W)^3\geq 3$. This is not possible in fact on these surfaces we have rk$(W)$=20 which implies rk$(W^{\vee}/W)^3\leq 2$. This means that the three classes fit together giving a 3-divisible class in NS($O_\mathcal{L}^{(8,2)}$) and NS($O_\mathcal{L}^{(8,3)}$). Now with the help of a blow-up it is possible to resolve the singularities of the total space $\pitr/O\times T$ of the family in correspondence of the singular fibers $O_\mathcal{L}^{(8,2)}$ and $O_\mathcal{L}^{(8,3)}$, obtaining a smooth threefold there, and where the two singular surfaces are replaced by surfaces with only $A_i$-singularities. Then by using a result of Briskorn, cf. \cite{briskorn}, one can replace these two fibers by smooth fibers, the total space remains smooth. Hence we see that the previous class which we have on the special fibers also exists in the Neron-Severi group of the generic fiber $O_{\mathcal{L}}$.

\eprf

\begin{rem}
In the proof of Proposition \ref{picard}, 2, we need a special resolution of the total space which induces a resolution of two special fibers. To do this we make a simultaneous blow-up of the threefold singularities and of the surface singularities. In fact we are able to do this in the case that the singularities of the surface is an $A_i$ or a $D_n$, $n$ even, singularity. We do not know how to find such blow up if the fibers have $D_n$, $n$ odd, singularities, and $E_6$, $E_7$ or $E_8$ singularities. One method to find the necessaries resolutions also in these cases is to use the famous McKay-correspondence, between group representations and resolutions of singularities (cf. \cite{inak}). This way works and it is the result of a forthcoming paper \cite{remc} of M. Lehn and the author.
\end{rem}

\subsection{The special cases}
\begin{prop}\label{picards}
1. The Picard-lattice of the special surfaces in $T_{\mathcal{L}}$ and $O_{\mathcal{L}}$ is generated in all the cases but 
$O_{\mathcal{L}}^{(8,4)}$ 
by the curves of  Proposition \ref{TO} and by  the 2-divisible classes of Proposition \ref{picard}. In the case of $O_{\mathcal{L}^{(8,4)}}$ the class:
\begin{eqnarray*}
L_1+L_3+L_5+N_1+C+N_4+R_2+M_1
\end{eqnarray*}
is 2-divisible too. \\
2. In the case of $T_{\mathcal{M}}^{(6,1)}$ and   of $T_{\mathcal{M}}^{(6,2)}$ the class:
\begin{eqnarray*}
\bar{L}:=N_1-N_2+N_3-N_4+N_5-N_6+N_7-N_8+N_9-N_{10}+N_{11}-N_{12}
\end{eqnarray*}
is 3-divisible. Moreover in the case of $T_{\mathcal{M}}^{(6,2)}$ the classes:
\begin{eqnarray*}
N_1+C_1+N_4+N_5+C_2+N_8+M_1+M_2,\\
N_1+C_1+N_4+N_9+C_3+N_{12}+M_1+M_3\\
\end{eqnarray*}
are 2-divisible. These together with the $20$ curves of Section \ref{special}  span a 20-dimensional lattice.\\
3.  In the case of  $O_{\mathcal{M}}^{(8,1)}$ and  $O_{\mathcal{M}}^{(8,4)}$, the class:
\begin{eqnarray*}
\bar{L}:=M_1+M_2+M_3+M_4+R_1+R_3+R_1'+R_3'
\end{eqnarray*}
is 2-divisible and in the case of  $O_{\mathcal{M}}^{(8,4)}$ the class:
\begin{eqnarray*}
W:=R_1+2R_2+3R_3+R_1'+2R_2'+3R_3'+2N_1+2C_1+3M_1+M_2+2N_3+2C_2+3M_3+M_4
\end{eqnarray*}
is 4-divisible.\\
This together with the 20 curves of Section \ref{special} spans a 20-dimensional lattice.

The ranks and the discriminants of the Picard-lattices are:
\begin{eqnarray*}
\begin{array}{c|c|c|c|c|c|c }
&\multicolumn{4}{c|}{T_{\mathcal{L}}}&\multicolumn{2}{c}{T_{\mathcal{M}}}\\
\hline
&6,1&6,2&6,3&6,4&6,1&6,2\\
\hline
{\rm rank}&20&20&20&20&20&20\\
\hline
{\rm discriminant}&-3\cdot 5&-2^2\cdot 3\cdot 5&-2^2\cdot 3\cdot 5&-3\cdot 5&
-3\cdot 5&-2^2\cdot 3\cdot 5
\end{array}
\end{eqnarray*}
\begin{eqnarray*}
\begin{array}{c|c|c|c|c|c|c}
&\multicolumn{4}{c|}{O_{\mathcal{L}}}&\multicolumn{2}{c}{O_{\mathcal{M}}}\\
\hline
&8,1&8,2&8,3&8,4&8,1&8,4\\
\hline
{\rm rank}&20&20&20&20&20&20\\
\hline
{\rm discriminant}&-2^2\cdot 7&-2^2\cdot 3\cdot 7&-2^3\cdot 3\cdot 7&-2^2\cdot 7&-2^2\cdot 7&-2^4\cdot 7
\end{array}
\end{eqnarray*}
\end{prop}
\bprf(cf. also \cite{basa} Theorem 6.2) 
1. In all the special cases of the pencils $T_{\mathcal{L}}$   and $O_{\mathcal{L}}$ we do not have enough rational $(-2)$-curves satisfying the conditions of \cite{nikulin} and \cite{tan} to obtain more 2-divisible or 3-divisible classes.\\
2. The class $\bar{L}$ comes from the 3-divisible class $\mathcal{L}$ on $Y_{T\times T}$ hence it is 3-divisible too. One proves as in Proposition \ref{picard} that in the case of $T_{\mathcal{M}}^{(6,2)}$ we have two more divisible classes.\\
3. The class $\bar{L}$ comes from the class $\mathcal{L}$ on  $Y_{O\times O}$, hence it is 2-divisible on $O_{\mathcal{M}}^{(8,1)}$ and  $O_{\mathcal{M}}^{(8,4)}$. 


Consider the classes:
\begin{eqnarray*}
\begin{array}{l}
v_1:=R_1+2R_2+3R_3,\\
v_2:=R_1'+2R_2'+3R_3',\\
v_3:=2N_1+2C_1+3M_1+M_2,\\
v_4:=2N_3+2C_2+3M_3+M_4\\
\end{array}
\end{eqnarray*}
on $O_{\mathcal{M}}^{(8,4)}$. Since the rank of NS is 20  the classes $v_i/4$, $i=1,2,3,4$ cannot be independent in $NS^{\vee}/NS$, hence we must find coefficients $\alpha_i\in\Z_4$ with: 
\begin{eqnarray*}
W:=1/4(\alpha_1v_1+\alpha_2v_2+\alpha_3v_3+\alpha_4v_4)\in NS.
\end{eqnarray*}
Observe that then $2W\in NS$ too and so the class: 
\begin{eqnarray*}
1/2(\alpha_1(R_1+R_3)+\alpha_2(R_2'+R_3')+\alpha_3(M_1+M_2)+\alpha_4(M_3+M_4))
\end{eqnarray*}
must be in $NS$ too. By \cite{nikulin} this must contains eight or 16 curves, hence $\alpha_i=1~(mod~2\Z)$ and after perhaps interchanging $v_i$ with $-v_i$, we may assume $\alpha_i=1~(mod~4\Z)$. In this way we find the class of the statement.   Finally by the same argumentation of 1 and  2 we cannot have more 3-divisible classes or  2-divisible classes.
\eprf
\section{More cyclic coverings}
\subsection{The general case}
As in Section \ref{cyclic} we can construct the 3-cyclic covering of $T_{\mathcal{L}}$, $T_{\mathcal{M}}$ by using the 3-divisible classes $\bar{L'}$ and $\bar{L}$ and the 2-cyclic coverings of $O_{\mathcal{L}}$, $O_{\mathcal{M}}$ by using the 2-divisible divisors $\bar{L'}$, $\bar{L}$ .  As before we consider first the non-special surfaces. Call the coverings $T_{\bar{L'}}$, $T_{\bar{L}}$, resp. $O_{\bar{L'}}$, $O_{\bar{L}}$. Then:
\begin{prop}
The surfaces  $T_{\bar{L'}}$ and $T_{\bar{L}}$ are isomorphic and the surfaces $O_{\bar{L'}}$ and $O_{\bar{L}}$  too. These are K3-surfaces which contains 15 rational curves coming from the 19 curves of $T_{\mathcal{L}}$, and 18 rational curves coming from the 19 curves of $O_{\mathcal{L}}$, moreover are isomorphic with  $Y_{ V\times V}$ and with  $Y_{ T\times T}$.
\end{prop}
\bprf
The irreducible curves contained on the divisor $\bar{L'}$ disappear on 
$T_{\bar{L'}}$ and the same happens on the surface $O_{\bar{L'}}$ . On these surfaces we have the configuration  of 15, resp. 18 $(-2)$-rational curves given in Figure 8.
\begin{figure}[t]

\input{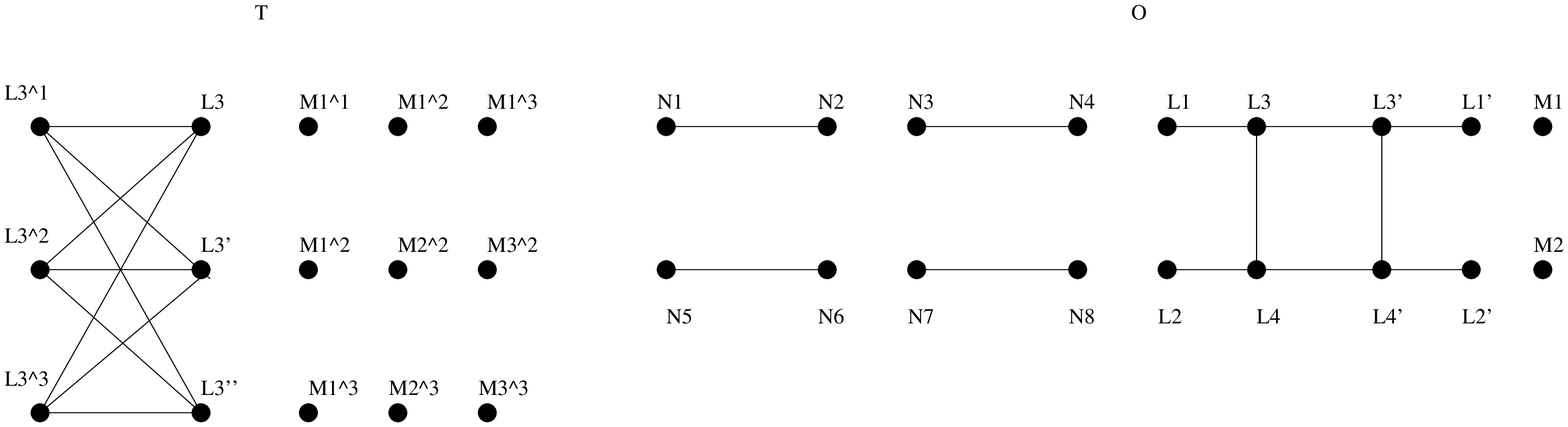}

\begin{center}
\small{Fig. 8}
\end{center}

\end{figure}

\noindent
Similarly to Proposition \ref{TO}  we write down the discriminants:
\begin{eqnarray*}
\begin{array}{c|ccc|c}
               &d(L)                &d(M)      &d(N)  &d\\
\hline
T_{\bar{L}'}   &-2^4\cdot 5         &-2^9      &      &2^{13}\cdot 5\\
O_{\bar{L}'}   &-3^2\cdot 7         &2^2      &3^4    &-2^2\cdot 3^6\cdot 7\\
\end{array}
\end{eqnarray*}

The surfaces are K3 and are isomorphic with $Y_{ V\times V}$ and with $Y_{ T\times T}$. Namely observe that $V\times V$ is a subgroup of index 3 of the groups $T\times V$ and of $(TT)'$, hence we have a diagram as those given in Figure 2. Similarly $T\times T$ is an index 2 subgroup of $O\times T$ and of $(OO)''$ so as before we have the diagram of Figure 2.
\eprf
\subsection{The special cases}\label{mcspecial}
We construct 3-cyclic  coverings of the special $T_{\mathcal{L}}$ and $T_{\mathcal{M}}$ by using the 3-divisible classes $\bar{L'}$ and $\bar{L}$ and the 2-cyclic covering of the special $O_{\mathcal{L}}$ and  $O_{\mathcal{M}}$ by using the 2-divisible classes $\bar{L'}$, $\bar{L}$ (in the same way as in Section \ref{cyclic}). Other than the curves coming from the base locus of $X^6_{\lambda}$ and $X^8_{\lambda}$ we have the configurations of rational curves given in Figure 9. Observe that as in the general case the surfaces $T_{\bar{L}'}$ are isomorphic to the surfaces $T_{\bar{L}}$ and the surfaces $O_{\bar{L}'}$ are isomorphic to the $O_{\bar{L}}$.
\begin{figure}[t]

\input{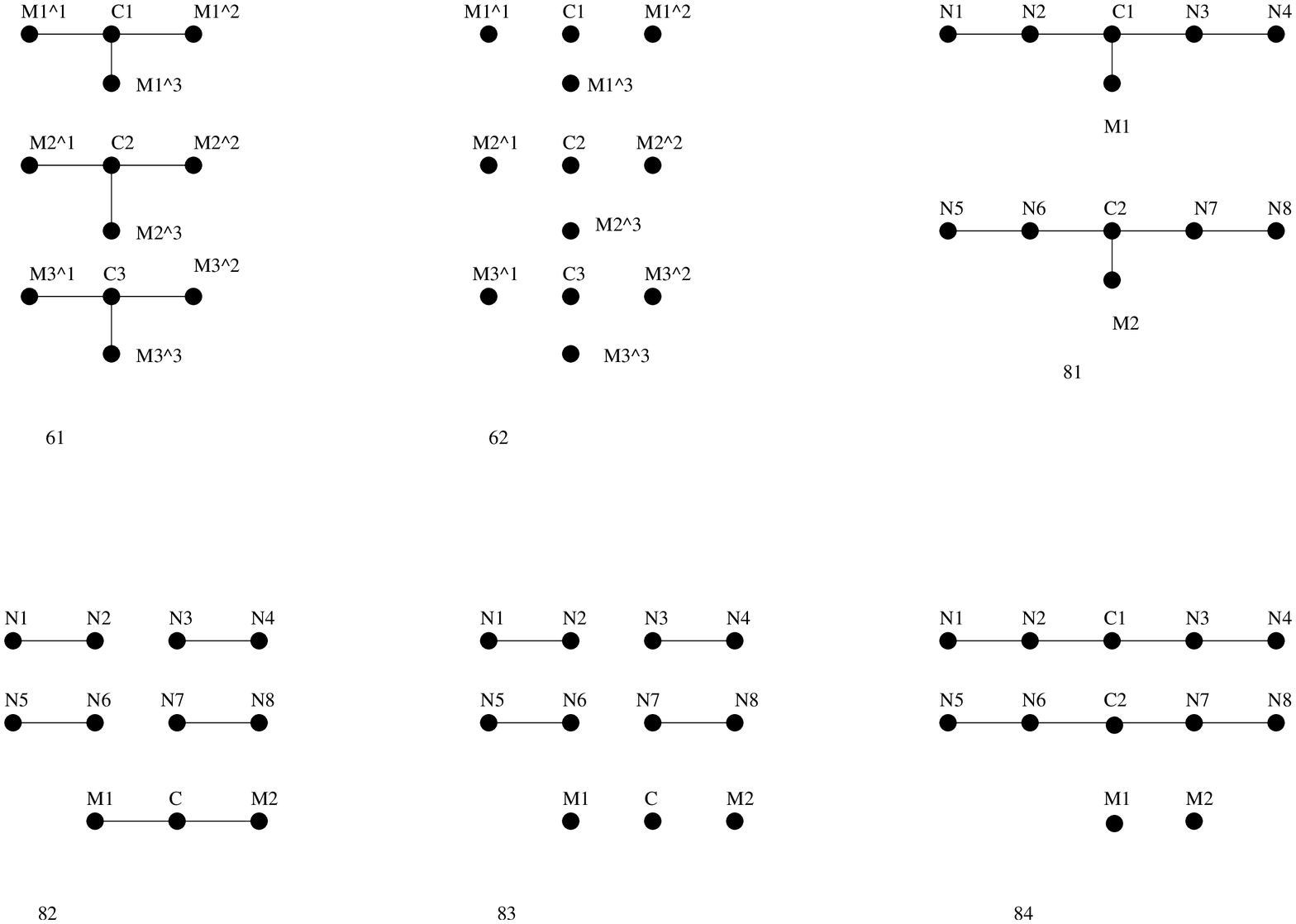}

\begin{center}
\small{Fig. 9}
\end{center}

\end{figure}

\noindent
The surfaces which we obtain are exactly the special K3-surfaces in the families  $Y_{\lambda, V\times V}$ and $Y_{\lambda, T\times T}$. The number of independent rational curves in the Neron-Severi group is:
\begin{eqnarray*}
\begin{array}{r|r|r|r|r|r|r|r|r}
\bar{L} (\bar{L'})&6,1&6,2&6,3&6,4&8,1&8,2&8,3&8,4\\
\hline
&18&18&18&18&20&19&19&20
\end{array}
\end{eqnarray*}
This gives a lower bound for the Picard-number, the discriminants of the lattices generated by these curves are given in Section \ref{rational}.

\subsection{The Picard-Lattices}
In this Section we identify the Picard-lattices of the surfaces $O_{\bar{L}'}^{(8,1)}$ and $O_{\bar{L}'}^{(8,4)}$ which have $\rho=20$. The proof of the following Proposition is left to the reader since it is very similar to the proof of Propositions \ref{picard} and \ref{picards}.

\begin{prop}
The classes:
\begin{eqnarray*}
k_1'=L_2-L_4+L_3'-L_1'+N_1-N_2+N_3-N_4+N_5-N_6+N_7-N_8,\\
k_1''=L_2'-L_4'+L_3-L_1+N_1-N_2-N_3+N_4+N_5-N_6-N_7+N_8\\
\end{eqnarray*}
are 3-divisible in $NS(O^{(8,1)}_{\bar{L}'})$, $NS(O^{(8,4)}_{\bar{L}'})$. Moreover in the case of $O_{\bar{L}'}^{(8,4)}$ the class:
\begin{eqnarray*}
\kappa=N_1+C_1+N_4+N_5+C_2+N_8+M_1+M_2\\
\end{eqnarray*}
is 2-divisible. These classes together with the $(-2)$-curves of Section \ref{mcspecial} span lattices of the following ranks and discriminants:
\begin{eqnarray*}
\begin{array}{r|r|r}
&8,1&8,4\\
\hline
{\rm rank}&20&20\\
\hline
{\rm discriminant}&-7&-2^2\cdot 7\\
\end{array}
\end{eqnarray*}
\end{prop}

\begin{rem}
Observe that the surfaces $Y_{\lambda, V\times V}$ and $Y_{\lambda, T\times T}$ are for each $\lambda\in\piu$ the fiber products
\begin{eqnarray*}
\begin{array}{c}
Y_{\lambda, T\times V}\times Y_{\lambda, (TT)'}~\mbox{over}~Y^6_{\lambda, T\times T}\\
Y_{\lambda, O\times T}\times Y_{\lambda, (OO)''} ~\mbox{over}~Y_{\lambda, O\times O},\\
\end{array}
\end{eqnarray*}
(to avoid confusion we denote here by  $Y^6_{T\times T}$ the minimal resolution of $X^6_{\lambda}/ {T\times T}$). 

\end{rem}
\section{Final remarks}
1. In the Sections 8 and 9 we identify explicitly the Picard-lattice of some K3-surfaces. It is our next aim to compute the transcendental lattices orthogonal to the Picard-lattices to classify the K3-surfaces. In particular by a result of Shioda and Inose, cf. \cite{si}, K3-surfaces with $\rho=20$ are classified by means of their transcendental lattice.\\
2. By a result of Morrison, cf. \cite{mo}, each K3-surface with $\rho=19$ or $20$ admits a so called Shioda-Inose structure. This means that there is a Nikulin-involution, an involution with eight isolated fix-points and the quotient is birational to a Kummer-surface. It would be desirable to have an explicit description of this structure for our surfaces.\\
3. The quotients 3-folds $\pitr/G$, $G=T\times T, O\times O, I\times I$ or a subgroup described in  Section  \ref{normal} seem to be Fano. We did not describe these spaces here. It would be interesting to have a global resolution of these spaces and to see our $K3$-surfaces as pencils on the  smooth 3-folds. A first result which is helpful in finding such resolutions is contained in \cite{remc}.


\vspace*{1.0cm}

{\footnotesize FB 17 Mathematik und Informatik\\
Universit\"at Mainz\\
Staudingerweg 9\\
55099 Mainz\\
Germany\\
e-mail:sarti@mathematik.uni-mainz.de.}

\end{document}